\documentclass[12pt,a4paper]{amsart}
\usepackage{amssymb,amsmath}

\headsep=1cm \numberwithin{equation}{section}
\newtheorem{theorem}{Theorem}[section]
\newtheorem{lemma}[theorem]{Lemma}
\newtheorem{proposition}[theorem]{Proposition}
\newtheorem{corollary}[theorem]{Corollary}

\theoremstyle{definition}
\newtheorem{definition}[theorem]{Definition}
\newtheorem{definition and remark}[theorem]{Definition and Remark}
\newtheorem{remark}[theorem]{Remark}
\newtheorem{example}[theorem]{Example}

\newcommand{\Ass}{\operatorname{Ass}}

\newcommand{\id}{\operatorname{id}}
\newcommand{\grade}{\operatorname{grade}}

\newcommand{\Z}{\operatorname{Z}}
\newcommand{\Spec}{\operatorname{Spec}}

\newcommand{\Ht}{\operatorname{ht}}

\newcommand{\Dp}{\operatorname{dp}}

\newcommand{\injdim}{\operatorname{injdim}}
\newcommand{\Ext}{\operatorname{Ext}}
\newcommand{\D}{\operatorname{D}}
\newcommand{\Supp}{\operatorname{Supp}}
\newcommand{\supp}{\operatorname{supp}}

\newcommand{\Hom}{\operatorname{Hom}}

\newcommand{\Ann}{\operatorname{Ann}}

\newcommand{\depth}{\operatorname{depth}}

\newcommand{\coker}{\operatorname{coker}}

\newcommand{\COH}{\operatorname{H}}

\newcommand{\lo}{\longrightarrow}
\newcommand{\fm}{\frak{m}}
\newcommand{\fp}{\frak{p}}
\newcommand{\fq}{\frak{q}}
\newcommand{\fa}{\frak{a}}

\newcommand{\fx}{\bf{x}}

\def\mapdown#1{\Big\downarrow\rlap
{$\vcenter{\hbox{$\scriptstyle#1$}}$}}

\begin{document}
\author[Esmkhani and Tousi]{Mohammad Ali Esmkhani and Massoud Tousi}
\title[Duality for a Cohen-Macaulay local ring]
{Duality for a Cohen-Macaulay local ring}

\address{M.A. Esmkhani, Department of Mathematics, Zanjan
University, P.O. Box 45195-313, Zanjan, Iran.}
\email{esmkhani@znu.ac.ir}

\address{M. Tousi, Department of Mathematics, Shahid Beheshti University,
Tehran, Iran-and-Institute for Studies in Theoretical Physics and
Mathematics, P.O. Box 19395-5746, Tehran, Iran.}
\email{mtousi@ipm.ir}

\subjclass[2000]{Primary 13C14, 13J10; Secondary 13D05, 13C13.}

\keywords{Canonical module, balanced big Cohen-Macaulay module,
Cohen-Macaulay local ring, complete module, Cousin complex.}
\thanks{This research was in part supported by a grant from
IPM (No. 87130214)}

\begin{abstract}
Let $(R,\fm)$ be a Cohen-Macaulay local ring. If $R$ has a
canonical module, then there are some interesting results about
duality for this situation. In this paper, we show that one can
indeed obtain similar these results in the case $R$ has not a
canonical module. Also, we give some characterizations of complete
big Cohen-Macaulay modules of finite injective dimension and by
using it some characterizations of Gorenstein modules over the
$\fm$-adic completion of $R$ are obtained.
\end{abstract}

\maketitle

\section{Introduction}

Throughout the paper, $R$ will denote a commutative Noetherian
local ring with nonzero identity and maximal ideal $\fm$. Unless
stated otherwise, the notation will be the same as used in [{\bf
2}]. When discussing the completion of $R$ or a module over $R$,
we will mean the $\fm$-adic completion.

The canonical module of a Cohen-Macaulay ring gives a useful
technical tool in commutative algebra, algebraic geometry, and
singularity theory via various duality theorems. Grothendieck
defined the notion of a module of dualizing differentials for a
complete local ring and proved some results on local duality for
complete local rings and some properties of a module of dualizing
differentials, see [{\bf 7}]. In [{\bf 8}], the theory of
canonical modules for a general Noetherian local ring $R$ was
studied. We recall the definition of the canonical modules of a
Noetherian local ring $(R,\fm)$ of dimension $n$. A finitely
generated $R$-module $C$ is called the canonical module of $R$ if
$C\otimes_R\hat{R}\cong \Hom_R(H_{\fm}^ n(R),E_R(R/\fm))$, see
[{\bf 8}]. There are many unknowns about canonical modules over
general local rings, existence theorems are among those. In [{\bf
10}], Reiten (and independently Foxby [{\bf 5}]) proved that a
Cohen-Macaulay local ring has a canonical module if and only if it
is a homomorphic image of a Gorenstein local ring. Therefore,
there exists a Cohen-Macaulay local ring which does not have
canonical modules. It is known that any two canonical modules of
$R$ are isomorphic. Also, it is known that over a Cohen-Macaulay
local ring $R$, a maximal Cohen-Macaulay module $C$ of type 1 and
of finite injective dimension is a canonical module of $R$.

Section 2 develops a variation of Grothendieck duality, namely
duality over a Cohen-Macaulay local ring. Let $R$ be a
Cohen-Macaulay local ring and $K$ the canonical module of the
completion of $R$. The $R$-module $K$ is not finitely generated in
general. Many facts that are valid for finitely generated modules
are not valid for $K$. But our achievement in this paper is that,
$K$ as an $R$-module has many properties similar to canonical
modules. We do this by applying some results of H.B. Foxby [{\bf
6}], A.M. Simon [{\bf 18}] and R.Y. Sharp [{\bf 15}]. Simon [{\bf
18}] has established analogies between complete modules and
finitely generated ones. The purpose of [{\bf 15}] is to show that
balanced big Cohen-Macaulay modules (over an arbitrary local ring
$R$) have many of the properties of maximal Cohen-Macaulay
modules. Foxby [{\bf 6}] has established some properties of the
depth of finitely
generated modules for a non-finitely generated modules. \\ \\
One of the main result of this paper is\\ \\
{\bf Theorem A}  {\it (See Theorem 2.9) (Duality) Let $(R,\fm)$ be
a Cohen-Macaulay local ring of dimension $d$, $K$ the canonical
module of the completion of $R$, and $\D$ the functor
$\Hom_R(-,K)$. Then

i) $D$ takes balanced big Cohen-Macaulay $R$-module to balanced
big Cohen-Macaulay $R$-module.

ii) $D$ takes exact sequence of balanced big Cohen-Macaulay
$R$-modules to exact sequence.

iii) There exists a natural transformation $\tau
:-\otimes_R\hat{R}\lo \Hom_R(\Hom_R(-,K),K)$ such that $\tau _M$
sending $m\otimes \hat{r}\in M\otimes_R\hat{R}$ to the map
$\alpha\rightarrowtail \hat{r}\alpha(m)$ for $\alpha\in
\Hom_R(M,K)$ is an isomorphism when $M$ is a maximal
Cohen-Macaulay $R$-module.}\\

The purpose of Section 3 and 4 is to present some
characterizations of $K$. We do this by characterizing the
Gorenstein modules over $\hat{R}$. Gorenstein $R$-modules are
those finitely generated $R$-modules for which the Cousin complex
for $M$ with respect to $M$-height filtration provides a minimal
injective resolution (see [{\bf 13}]). R.Y. Sharp obtained various
characterizations and properties of Gorenstein modules in [{\bf
13}, Theorem 3.1]. Note that, it is shown that in [{\bf 14}] that
any Gorenstein $\hat{R}$-module is a finite sum of copies of $K$.
First, we obtain some characterizations for complete big
Cohen-Macaulay modules (see Proposition 3.2). For complete big
Cohen-Macaulay modules of finite injective dimension a better
characterization is
provided:\\ \\
{\bf Theorem B} {\it (See Theorem 3.4) Let $(R,\fm)$ be a
Noetherian local ring and $M$ a nonzero complete $R$-module. Then
the following conditions are equivalent:\\
i) $\Ht_M(\fm)=\dim(R)$ and for every $\fp \in \supp_R(M)$,
$\mu_R^ i(\fp, M)\neq 0$ if and only if $i\neq
\Ht_{M_{\fp}}({\fp}R_{\fp})$.\\
ii) $\mu_R^ i(\fm, M)=0$ if and only if $i\neq \dim(R)$.\\
iii) $M$ is a big Cohen-Macaulay $R$-module with respect to a
system of
parameters for $R$ and $M$ has finite injective dimension.\\
iv) $\depth_R(M)=\dim(R)$ and $M$ has finite injective
dimension.\\
v) $\injdim_R(M)=\depth_R(M)=\dim_R(M)=\dim(R)=\depth(R)$.}\\ \\
Finally, we establish the following theorem\\ \\
{\bf Theorem C} {\it (See Theorem 4.2) Let $(R,\fm)$ be a
Noetherian local ring of dimension $d$ and $M$ a nonzero
$R$-module. Then the following conditions are
equivalent:\\
i) $M$ is complete and big Cohen-Macaulay with respect to a system
of parameters $x_1,\ldots,x_d$ such that $\injdim_RM$ and $\mu_R^
d(\fm,M)$ are finite.\\
ii) $M$ is a Gorenstein $\hat{R}$-module.\\
iii) $M$ is complete and for all Cohen-Macaulay $R$-module $L$ of
dimension $t$, we have

\ \ \ a) $\Ext_R^ i(L,M)=0$ for all $i\neq d-t$, and

\ \ \ b) $\Ext_R^{d-t}(\Ext_R^{d-t}(L,M),M)\cong L\otimes_RF$,
where $F$ is a finitely generated free $\hat{R}$-module.}

\section{The canonical modules of the completion of a Cohen-Macaulay local ring}

We begin this section, by recalling the definitions of the notion
of depth, dimension of an $R$-module $M$, and balanced big
Cohen-Macaulay modules (recall that $M$ is not assumed to be
finitely generated).

\begin{definition and remark} Let $(R,\fm)$ be a Noetherian local
ring.

i) If $M$ is an $R$-module, we define the dimension of $M$ by
$$\dim_R(M)=\dim_R(R/{{\Ann_R(M)}}).$$ By the height $\Ht_M({\fm})$
of an $R$-module $M$ we mean the dimension of the support
$\Supp_R(M)$ of $M$.

ii) In [{\bf 6}, Section 1], Foxby has established the concept of
depth to any $R$-module $M$ by means of the formula
$$\depth_R(M)=\inf\{i\geq 0\mid \Ext_R^i(R/\fm,M)\neq 0\}.$$ For a
nonzero $M$ which is not finitely generated, it is possible for
this depth to be $\infty$. If the $R$-module $M$ satisfies $\fm
M\neq M$, then $\depth_R(M)\leq \Ht_M({\fm})< \infty$ (see [{\bf
6}, Corollary 2.2]).

iii) (See [{\bf 19}, Definition 5.3.1]) Let $M$ be an $R$-module
with $\fm M\neq M$. Then $\Dp_R(M)$ is defined to be the supremum
of the lengths of all $M$-sequence contained in $\fm$. Since $R$
is Noetherian, $\Dp_R(M)< \infty$. It is well known that if $M$ is
nonzero finitely generated, then $\Dp_R(M)=\depth_R(M)$.

iv) Let $M$ be an $R$-module with $\fm M\neq M$. Then, by ii) and
[{\bf 19}, Proposition 5.3.7 (ii)], $$\Dp_R(M)\leq \depth_R(M)\leq
\Ht_M({\fm})\leq \dim_R(M)\leq \dim(R).$$

v) Let $a_1,\ldots,a_n$ be a system of parameters (s.o.p.) for
$R$. A (not necessarily finitely generated) $R$-module $M$ is said
to be a big Cohen-Macaulay $R$-module with respect to
$a_1,\ldots,a_n$ if $a_1,\ldots,a_n$ is an $M$-sequence, that is
if $M\neq (a_1,\ldots,a_n)M$ and, for each $i=1,\ldots,n$,
$((a_1,\ldots,a_{i-1})M: a_{i})=(a_1,\ldots,a_{i-1})M$. It will be
convenient to use abbreviation b.C.M. for big Cohen-Macaulay.

vi) Let $a_1,\ldots, a_n$ be a system of parameters for $R$ and
let $M$ be b.C.M. $R$-module with respect to $a_1,\ldots, a_n$.
Then there exists $t\in \mathrm{N}$ such that ${\fm}^ t\subseteq
(a_1,\ldots, a_n)R$ and $a_1,\ldots, a_n$ is an $M$-sequence. So,
$\fm M\neq M$. Hence, by iv), we have
$$\Dp_R(M)=\depth_R(M)=\Ht_M({\fm})=\dim_R(M)=\dim(R).$$

vii) An $R$-module $M$ is called balanced big Cohen-Macaulay if
every system of parameter for $R$ is an $M$-sequence (see [{\bf
15}, Definition 1.4]). It will be convenient to use abbreviation
b.b.C.M. for balanced big Cohen-Macaulay. In [{\bf 15,16}], it is
shown that b.b.C.M. modules verify some of the properties of
finitely generated Cohen-Macaulay modules.

\end{definition and remark}

\begin{lemma} Let $(R,\fm)$ be a Cohen-Macaulay local ring and $K$ the canonical module of $\hat{R}$. Then

i) $K$ is complete with respect to $\fm$-adic topology.

ii) $K$ is a b.b.C.M. $R$-module.

iii) $\injdim_RK< \infty$.

\end{lemma}

{\bf Proof.} i) The natural $R$-homomorphism $\theta :K\lo
{\hat{K}}$ is an isomorphism. This follows from the fact that
$\fm^{n}K=(\fm{\hat{R}})^{n}K$ for all $n\geq 0$.

ii) Let $x_1,...,x_d$ be a system of parameters for $R$ and let
$\psi:R\lo {\hat{R}}$ be the natural ring homomorphism. Then
$\psi(x_1),...,\psi(x_d)$ is an ${\hat{R}}$-sequence, because $R$
is Cohen-Macaulay. On the other hand, $K$ is a maximal
Cohen-Macaulay ${\hat{R}}$-module, and so $x_1,...,x_d$ is a
$K$-sequence.

iii) This follows from the fact that any injective
$\hat{R}$-module is injective as an $R$-module. Note that the
injective dimension ok $K$ as an ${\hat{R}}$-module is finite.
$\Box$

Simon [{\bf 18}] has shown that complete modules behave similar to
finitely generated modules in many respects. We summarize some
important properties of complete modules from [{\bf 18}] in the
following remark.
\begin{remark} Let $(R,\fm)$ be a  Noetherian local ring, and let
$M$ and $N$ be two complete $R$-modules. Then

i) If $f: N\lo M$ is a morphism between $R$-modules such that
$M=f(N)+\fm M$, then $f$ is surjective.

ii) For each $R$-module $X$, if $\Ext_R^i(X,M)\neq 0$, then
$\Ext_R^i(X,M)\neq \fm\Ext_R^i(X,M)$.

iii) $\injdim_R(M)=\sup\{i\mid \Ext_R^i(R/\fm,M)\neq 0\}$.

iv) If $M$ has finite injective dimension, then
$\injdim_R(M)=\depth R$.

v) If $R$ has a big Cohen-Macaulay module $C$ and a complete
module $M$ of finite injective dimension, then the ring $R$ is
Cohen-Macaulay.
\end{remark}

To prove Theorem 2.9, we need the following five results.

\begin{lemma} Let $(R,\fm)$ be a Noetherian local ring and
$M$ an $R$-module of dimension $t\geq 1$ with $\fm M\neq M$ ($M$
not necessarily finitely generated). If there exists an
$M$-sequence $x_1,...,x_t$, then
$$\dim_RM=\dim_R(M/{x_{1}M})+1.$$
\end{lemma}

{\bf Proof.} It is clear that $x_{1}\in R \setminus \Z_R(M)$
implies that $x_{1}\in R\setminus \Z_R(R/{\Ann_RM})$. Hence, we
have
\begin{equation*}
\begin{split}
\depth_R(M/{x_{1}M})\leq \dim_R(R/{(\Ann_R(M/{x_{1}M}))})&\leq
\dim_R(R/{(\Ann_RM+x_{1}R)})\\&=
\dim_R({(R/{\Ann_RM})}/{x_{1}{(R/{\Ann_RM})}})\\&=\dim_R(R/{\Ann_RM})-1\\&=t-1,
\end{split}
\end{equation*}
where the first inequality follows from [{\bf 6}, Corollary 2.2].
On the other hand, by [{\bf 2}, Proposition 1.2.3 and Lemma
1.2.4], $\Ext_R^i(R/{\fm},M/{x_{1}M})=0$ for all $0\leq i\leq
{t-2}$. So $\depth_R(M/{x_{1}M})\geq t-1$. Thus, the assertion
follows. $\Box$

\begin{proposition} Let $(R,\fm)$ be a Noetherian local ring.

i) Let $C$ be a complete $R$-module of injective dimension $s$. If
$M$ is an $R$-module and $x_1,x_2,...,x_l$ an $M$-sequence, then
$\Ext_R^j(M,C)=0$ for $j>s-l$.

ii) Let $(R,\fm)$ be a Cohen-Macaulay local ring of dimension $d$
and $C$ a b.b.C.M. $R$-module. If $M$ is an $R$-module of
dimension $t$, then  $\Ext_R^j(M,C)=0$ for all $j<d-t$.

iii)  Let $(R,\fm)$ be a Cohen-Macaulay local ring of dimension
$d$, $C$ a complete $R$-module of finite injective dimension which
is b.b.C.M. and $M$ an $R$-module of dimension $t\geq 0$. Let
$x_1,x_2,...,x_t$ be an $M$-sequence. If $\fm M\neq M$, then
$\Ext_R^{d-t}(M,C)\neq 0$, $x_1,x_2,...,x_t$ is an
$\Ext_R^{d-t}(M,C)$-sequence and the dimension of
$\Ext_R^{d-t}(M,C)$ is equal to $t$.
\end{proposition}

{\bf Proof.} i) The proof is same as [{\bf 2}, 3.3.3(b)], by
replacing Nakayama's lemma by Remark 2.3 ii).

ii) We set $\fa=\Ann_RM$. Since $R$ is Cohen-Macaulay, we have
$$\Ht_R({\fa})=\dim R-\dim{R/\fa}=d-t.$$ Hence $\grade_R({\fa})=d-t$, and
so there exist $x_1,x_2,...,x_{d-t}$ in $\fa$ which is an
$R$-sequence. Consequently, $x_1,x_2,...,x_{d-t}$  is a
$C$-sequence in $\fa$. Therefore, by [{\bf 2}, Proposition 1.2.3
and Lemma 1.2.4], we have $\Ext_R^j(M,C)=0$ for all $j<d-t$.

iii) We set $\Ext_R^{d-t}(M,C)=N$ and use induction on $t$. Let
$t=0$. There exists a system of parameters $y_1,y_2,...,y_d$ of
$R$ in $\Ann_R(M)$. So, by [{\bf 2}, Lemma 1.2.4],
$$\Ext_R^d(M,C)\cong \Hom_R(M,C/{(y_1,y_2,...,y_d)C}).$$
Since $C/{(y_1,y_2,...,y_d)C}\neq 0$ and
$(y_1,y_2,...,y_d)R\subseteq \Ann_R(C/{(y_1,y_2,...,y_d)C})$, we
have $\fm \in \Ass_R(C/{(y_1,y_2,...,y_d)C})$. Hence, there is an
$R$-monomorphism $\varphi: R/{\fm}\lo C/{(y_1,y_2,...,y_d)C}$. By
viewing $M/{\fm M}$ as a nonzero vector space over $R/{\fm}$, we
deduce that there is an $R$-epimorphism $\psi: M/{\fm M}\lo
R/\fm$. Hence $\varphi \psi \pi:M\lo C/{(y_1,y_2,...,y_d)C}$ is a
nonzero homomorphism in which $\pi:M\lo M/\fm M$ is the natural
$R$-epimorphism. So, $\Ext_R^d(M,C)\neq 0$. Also,
$\Ann_RM\subseteq\Ann_RN$. Therefore, $\dim_RN=0$.

Now suppose, inductively that $t>0$.  It is clear that $\fm
(M/x_1M)\neq (M/x_1M)$. Now by Lemma 2.4, $\dim_R(M/x_1M)=t-1$. In
view of i), ii) and the fact that
$\injdim_R(C)=\depth(R)=\dim(R)$, the exact sequence $0\lo
M\overset{x_1}\lo M\lo M/{x_1M}\lo 0$ induces the exact sequence
$$0\lo N\overset{x_1}\lo N\lo \Ext_R^{d-t+1}(M/x_1M,C)\lo 0.
$$
So that, by inductive hypothesis, $\Ext_R^{d-t+1}(M/{x_1M},C)\neq
0$, $x_1,x_2,...,x_t$ is an $N$-sequence and the dimension of
$\Ext_R^{d-t+1}(M/x_1M,C)$ is $t-1$. Also, since $N\neq 0$, it
follows from Remark 2.3 ii) that $\fm N\neq N$. So, by Lemma 2.4,
$\dim_RN=t$, because $\dim_R(N/{x_1N})=t-1$. $\Box$ \\

In the following we state a result without proof, since the proof
is same as that of [{\bf 3}, Proposition 21.12].

\begin{lemma} Let $(R,\fm)$ be a Noetherian local ring of
dimension $d$, and $C$ a b.b.C.M. $R$-module. If $M$ is an
$R$-module and $\Ext_R^1(M,C)=0$, then for any $R$-sequence $x$ we
have
$$\Hom_R(M,C)/x\Hom_R(M,C) \cong \Hom_{R/xR}(M/xM,C/xC)$$
by the homomorphism taking the class of map $\varphi:M\lo C$ to
the map $M/xM \lo C/xC$ induced by $\varphi$.
\end{lemma}

\begin{lemma} Let $(R,\fm)$ be a Noetherian local ring, and
let $M$ and $N$ be two complete $R$-modules. Suppose that $x\in
\fm$ is a nonzero-divisor on $M$, $\varphi:N\rightarrow M$ an
$R$-homomorphism and $\overline{\varphi}:N/xN\rightarrow
 M/xM$ is the $R$-homomorphism induced by $\varphi$. Then\\
i) If $\overline{\varphi}$ is an epimorphism, then $\varphi$ is
an epimorphism.\\
ii) If $\overline{\varphi}$ is a monomorphism, then $\varphi$ is a
monomorphism.
\end{lemma}

{\bf Proof.} i) If $\overline{\varphi}$ is an epimorphism, then
$\varphi(N)+xM=M$. The assertion follows from Remark 2.3 i).\\
ii) Suppose that $\overline{\varphi}$ is a monomorphism, and let
$K=\ker \varphi$. Then $xK=K$, by same argument to the proof of
[{\bf 3}, Proposition 21.13]. So $K=\bigcap_{n=0}^\infty x^n
K\subseteq \bigcap_{n=0}^\infty x^n N$. But, since $N$ is a
complete $R$-module, we have $\bigcap_{n=0}^\infty x^n N=0$, and
so $K=0$. $\Box$

\begin{proposition} Let $(R,\fm)$ be a Cohen-Macaulay local ring
of dimension $d$ and $K$ the canonical module of the completion of
$R$. The following assertions hold.

i) $\Hom_R(H_\fm^d(R),E_R(R/\fm))\cong K$.

ii) $H_\fm^d(K)\cong E_R(R/\fm)$.

iii) There is an $R$-isomorphism $\psi :\hat{R}\lo \Hom_R(K,K)$
which is such that $\psi(\hat{r})=\hat{r}\id_K$, for any
$\hat{r}\in \hat{R}$.
\end{proposition}

{\bf Proof.} i) By using the Flat Base Change Theorem for local
cohomology, we obtain
\begin{equation*}
\begin{split}
 K&\cong
\Hom_{\hat{R}}(H_{\fm\hat{R}}^ d(\hat{R}
),E_{\hat{R}}(\hat{R}/\fm\hat{R}))\\&\cong
\Hom_{\hat{R}}(H_{\fm}^d(R)\otimes_R\hat{R},E_{\hat{R}}(\hat{R}/\fm\hat{R}))
\\&\cong \Hom_R(H_\fm^d(R),\Hom_{\hat{R}}(\hat{R},E_{\hat{R}}(\hat{R}/\fm\hat{R}))
\\&\cong \Hom_R(H_{\fm}^d(R),E_R(R/{\fm })).
\end{split}
\end{equation*}
\\
ii) In view of the Independence Theorem for local cohomology and
[{\bf 1}, Theorem 11.2.8], we have $H_\fm^d(K)\cong
H_{\fm\hat{R}}^d(K)\cong E_{\hat{R}}(\hat{R}/\fm\hat{R})\cong
E_R(R/\fm)$. \\ iii) We have
\begin{equation*}
\begin{split}
\Hom_R(K,K)&\cong \Hom_R(K,\Hom_R(H_\fm^d(R),E_R(R/\fm)))\\&\cong
\Hom_R(K\otimes_RH_\fm^d(R),E_R(R/\fm))\\&\cong
\Hom_R(H_\fm^d(K),E_R(R/\fm))\\&\cong
\Hom_R(E_R(R/\fm),E_R(R/\fm))\\&\cong \hat{R}.
\end{split}
\end{equation*}
Here the third isomorphism comes from [{\bf 1}, Ex.6.1.9]. So,
$\Hom_R(K,K)$ is complete.

Now, we prove this part by induction on $d=\dim R$. When $d=0$,
$R$ is complete, and so we have the result by [{\bf 2}, Theorem
3.3.4(d)]. Now suppose that $d> 0$ and the result has been proved
for the smaller values. Let $x\in \fm$ be a nonzero-divisor on
$R$. Then the module $K/xK$ is the canonical module of
$\hat{R}/{x\hat{R}}$. Also, there exists a natural isomorphism
$\hat{R}/{x\hat{R}}\cong \widehat{(R/xR)}$. Hence, by the
inductive hypothesis, we have an $R/xR$-isomorphism
$$\theta:\hat{R}/{x\hat{R}}\lo \Hom_{R/xR}(K/xK,K/xK)$$ which is
such that $\theta(\hat{r}+x\hat{R})=\hat{r}\id_{K/xK}$, for all
$\hat{r}\in \hat{R}$. It follows from Lemma 2.2, Remark 2.3 and
Proposition 2.5 i) that $\Ext^1_R(K,K)=0$. So, by Lemma 2.2 ii)
and Lemma 2.6, we have an $R$-isomorphism
$$f:\Hom_R(K,K)/x \Hom_R(K,K)\lo \Hom_{R/xR}(K/xK,K/xK)$$ by the homomorphism
taking the class $\varphi:K\lo K$ to the map $K/xK\lo K/xK$
induced by $\varphi$. Let $\overline{\psi}:\hat{R}/{x\hat{R}}\lo
\Hom_R(K,K)/x \Hom_R(K,K)$ be the homomorphism induced by
${\psi}$. It is easy to see that $\theta=f\overline{\psi}$. Hence
$\overline{\psi}$ is an $R$-isomorphism. The assertion follows
from Lemma 2.7. $\Box$
\\

Now, we are ready to prove one of the main results of this paper.

\begin{theorem} (Duality) Let $(R,\fm)$ be a Cohen-Macaulay local ring of dimension $d$, $K$
the canonical module of the completion of $R$, and $\D$ the
functor $\Hom_R(-,K)$. Then

i) $D$ takes b.b.C.M. $R$-module to b.b.C.M. $R$-module.

ii) $D$ takes exact sequence of b.b.C.M. $R$-modules to exact
sequence.

iii) There exists a natural transformation $\tau
:-\otimes_R\hat{R}\lo \Hom_R(\Hom_R(-,K),K)$ such that $\tau _M$
sending $m\otimes \hat{r}\in M\otimes_R\hat{R}$ to the map
$\alpha\rightarrowtail \hat{r}\alpha(m)$ for $\alpha\in
\Hom_R(M,K)$ is an isomorphism when $M$ is a maximal
Cohen-Macaulay $R$-module.
\end{theorem}

{\bf Proof.} i) Let $M$ be a b.b.C.M. $R$-module. Then, by Remark
and definition 2.1 vi), $\fm M\neq M$ and $\dim_R(M)=d$. So, the
claim follows from Lemma 2.2 and Proposition 2.5 iii).

ii) If $0\lo M'\lo M\lo M^{''}\lo 0$ is an exact sequence of
b.b.C.M. $R$-modules, then we get an exact sequence $$0\lo
\Hom_R(M^{''},K)\lo \Hom_R(M,K)\lo \Hom_R(M',K)\lo
\Ext_R^1(M^{''},K).$$ Using Lemma 2.2, Remark 2.3 iv) and
Proposition 2.5 i), we obtain $\Ext_R^1(M^{''},K)=0$. So, we have
the exact sequence $$0\lo D(M^{''})\lo D(M)\lo D(M')\lo 0.$$

iii) By Proposition 2.8 iii), there is an isomorphism $\psi
:\hat{R}\lo \Hom_R(K,K)$ which is such that
$\psi(\hat{r})=\hat{r}\id_K$ for any $\hat{r}\in \hat{R}$. Also,
we know that there is the natural transformation of functors
$\sigma :-\otimes_R\Hom_R(K,K)\lo \Hom_R(\Hom_R(-,K),K)$ such that
for each $R$-module $M$, $(\sigma_M(m\otimes f))(g)=f(g(m))$, for
$m\in M$, $f\in \Hom_R(K,K)$, and $g\in \Hom_R(M,K)$ (see [{\bf
1}, Lemma 10.2.16]). We define a natural transformation of
functors
$$\tau :-\otimes_R\hat{R}\lo \Hom_R(\Hom_R(-,K),K)$$ which is such
that for each $R$-module $M$, $\tau _M=\sigma _M(\id_M\otimes_R
\psi)$. Let $M$ be a maximal Cohen-Macaulay $R$-module and let
$F_1\overset{g}\lo F_0\overset{f}\lo M\lo 0$ be
an exact sequence where $F_0$ and $F_1$ are finitely generated free $R$-modules.\\
We consider the following cases:\\
Case 1. $\ker(f)\neq 0,\ker(g)\neq 0$. We consider two exact
sequences $$0\lo \ker(f)\overset{\rho}\lo F_0\overset{f}\lo M\lo
0,$$
$$0\lo \ker(g)\overset{\mu}\lo F_1\overset{g}\lo \ker(f)\lo 0.$$
It is easy to see that $\ker(f)$ and $\ker(g)$ are maximal
Cohen-Macaulay $R$-modules. By the above facts, i) and ii), we
obtain the following two exact sequences
$$0\lo \D(\D(\ker(f)))\lo \D(\D(F_0))\lo \D(\D(M))\lo 0$$ and
$$0\lo \D(\D(\ker(g)))\lo \D(\D(F_1))\lo \D(\D(\ker(f)))\lo 0.$$
Therefore, we have the exact sequence $$\D(\D(F_1))\lo
\D(\D(F_0))\lo \D(\D(M))\lo 0.$$ We obtain the following diagram

\begin{equation*}
\setcounter{MaxMatrixCols}{11}
\begin{matrix}
&F_1\otimes_R\hat{R} &{\lo} &F_0\otimes_R\hat{R} &{\lo}
&M\otimes_R\hat{R} &\lo &0
\\&\mapdown{\tau_{F_1}}& & \mapdown{\tau_{F_0}}
 & &\mapdown{\tau_M} & & & &
\\  &\D(\D(F_1)) &{\lo}
&\D(\D(F_0))&{\lo} &\D(\D(M)) &\lo &0.
\end{matrix}
\end{equation*}
Now, by using [{\bf 11}, Lemma 3.59], we get that $\sigma_{F_0}$
and $\sigma_{F_1}$ are isomorphism and therefore $\tau_{F_0}$ and
$\tau_{F_1}$ are also isomorphism. So,
$\tau_M$ is an isomorphism.\\
Case 2. Either $\ker(f)=0$ or $\ker(g)=0$. In view of the above
argument, it is obvious. $\Box$
\\

The following example shows that Theorem 2.9 iii) does not hold
for b.b.C.M. $R$-modules.

\begin{example} Let $(R,\fm)$ be a Gorenstein local ring
of dimension zero and $F$ a non-finitely generated free
$R$-module. Then the canonical module of $R$ is $R=E_R(R/\fm)$ and
$F$ is a b.b.C.M. $R$-module. Suppose that $\tau_F$ is an
isomorphism. This means that $F$ is Matlis reflexive and so by
[{\bf 4}, Proposition 1.3], we have an exact sequence $$0\lo L\lo
F\lo A\lo 0,$$ where $A$ is an Artinian $R$-module and $L$ a
finitely generated $R$-module. Therefore, $\mu_R^ 0(\fm,F)$ is
finite. This is a contradiction.
\end{example}

The following corollary is a generalization of Theorem 2.9.

\begin{corollary} Let $(R,\fm)$ be a Cohen-Macaulay local ring of dimension $d$ and
let $K$ be the canonical module of the completion of $R$. For all
integers $t=0,1,\ldots,d$ and all $R$-modules $M$ of dimension t
which is such that $\fm M\neq M$ and there exists an $M$-sequence
$x_1,\ldots,x_t$ we have

i) $\Ext_R^{d-t}(M,K)$ is an $R$-module of dimension $t$ and
$x_1,\ldots,x_t$ is an $\Ext_R^{d-t}(M,K)$-sequence,

ii) $\Ext_R^{i}(M,K)=0$ for all $i\neq d-t$, and

iii)if $M$ is finitely generated, then there exists an isomorphism
$M\otimes_R\hat{R}\lo \Ext_R^{d-t}(\Ext_R^{d-t}(M,K),K)$ which in
the case $d=t$ is just the natural isomorphism from
$M\otimes_R\hat{R}$ into the bidual of $M$ with respect to $K$.
\end{corollary}

{\bf Proof.} We obtain i) and ii) from Lemma 2.2 and Proposition
2.5.

iii) For $t=d$, the assertion follows from Theorem 2.9. By [{\bf
9}, page 140, Lemma 2(i)], we deduce that
$$\Ext_R^{d-t}(\Ext_R^{d-t}(M,K),K)\cong
\Hom_{R/{\fx R}}(\Hom_{R/\fx R}(M,K/{{\fx} K}),K/{{\fx} K})$$ for
an $R$-sequence ${\fx}$ of length $d-t$ which is contained in
$\Ann_R(M)$. Since $\dim_{R/{\fx} R}(M)=t$ and $\dim{R/{\fx}
R}=t$, so $M$ is a maximal Cohen-Macaulay $R/{\fx} R$-module. The
module $K/{\fx}K$ is the canonical module of $\hat{R}/{\fx}
\hat{R}$. Therefore, by Theorem 2.9 and above isomorphism, there
is a natural isomorphism
$$M\otimes_{R/{\fx} R}\widehat{R/{\fx} R}\lo
\Ext_R^{d-t}(\Ext_R^{d-t}(M,K),K).$$ In view of $\widehat{R/{\fx}
R}\cong {\hat{R}}/{\fx}{\hat{R}}\cong R/{\fx} R\otimes_R\hat{R}$,
we obtain the result. $\Box$

\section{Some characterizations for complete big Cohen-Macaulay modules of finite injective dimension}

We need to recall some definitions.

\begin{definition} Let $R$ be a Noetherian ring and let $M$ be an
$R$-module.

i) The small support, or little support, of $M$ denoted by
$\supp_R(M)$, is defined by
$$\supp_R(M)=\{\fp\in \Spec(R) \mid \depth_{R_{\fp}}(M_{\fp})<
\infty \}.$$ It is clear that $\supp_R(M)\subseteq \Supp_R(M)$; if
$M$ is finitely generated, then these two sets are equal, but in
general this need not be the case.

ii) A filtration of $\Spec(R)$ [{\bf 16}, 1.1] is a descending
sequence $\mathcal{F}=(F_i)_{i\geq 0}$ of subsets of $\Spec(R)$,
so that
$$F_0\supseteq F_1\supseteq \cdots \supseteq F_i\supseteq F_{i+1}\cdots,$$ with the property
that, for each $i\geq 0$, each member of $\partial
F_i=F_i\setminus F_{i+1}$ is a minimal member of $F_i$ with
respect to inclusion. We say that $\mathcal{F}$ admits $M$ if
$\Supp_R(M)\subseteq F_0$. Suppose $\mathcal{F}$ is a filtration
of $\Spec(R)$ that admits $M$. The Cousin complex $C(\mathcal{F}
,M)$ for $M$ with respect to $\mathcal{F}$  has the form
$$0\overset{d^{-2}}\lo M\overset{d^{-1}}\lo M^{0}\overset{d^{0}} \lo M^{1}\lo
\cdots \lo M^{n}\overset{d^{n}} \lo M^{n+1}\lo \cdots$$ with
$M^{n}=\oplus_{\fq\in \partial F_n}(\coker d^{n-2})_{\fq}$ for all
$n\geq 0$. The homomorphisms in this complex have the following
properties: for $m\in M$ and $\fq\in
\partial F_0$, the component of $d^{-1}(m)$ in $M_{\fq}$ is
$m/1$; for $n> 0$, $x\in M^{n-1}$ and $\fq\in \partial F_n$, the
component of $d^{n-1}(x)$ in $(\coker d^{n-2})_{\fq}$ is
$\overline{x}/1$, where $^- :M^{n-1}\lo \coker d^{n-2}$ is the
canonical epimorphism. The fact that such a complex can be
constructed is explained in [{\bf 16}, 1.3].

\end{definition}

In the following proposition, we obtain some characterizations for
complete big Cohen-Macaulay modules.

\begin{proposition} Let $(R,\fm)$ be a Noetherian local ring and $M$ a
nonzero complete $R$-module. Denote by
$\mathcal{D}(R)=(D_i)_{i\geq 0}$ (resp.
$\mathcal{H}(M)=(H_i)_{i\geq 0}$) the dimension filtration of \
$\Spec (R)$ (resp. $M$-height filtration of \ $\Spec (R)$) i.e.
$D_i=\{\ \fp \in \Spec (R)\mid \dim(R/{\fp})\leq \dim(R)-i\}$
(resp. $H_i=\{\ \fp \in \Supp_R(M)\mid
\Ht_{M_{\fp}}(\fp{R_{\fp}})\geq i\}$). Then the following conditions are equivalent:\\
i) $M$ is a b.C.M. $R$-module with respect to an s.o.p. for $R$.\\
ii) $\Dp_{R}(M)=\dim_R(M)=\Ht_M({\fm})=\depth_R(M)=\dim(R)$.\\
iii) $\depth_R(M)=\dim(R)$.\\
iv) $M$ is a b.b.C.M. $R$-module.\\
v) $C(\mathcal{D}(R), M)$ is exact.\\
vi) Both the complexes $C(\mathcal{D}(R), M)$ and
$C(\mathcal{H}(M), M)$ are isomorphic and $C(\mathcal{H}(M), M)$
is exact.\\
vii) $C(\mathcal{H}(M), M)$ is exact and $\Ht_M({\fm})=\dim(R)$.
\end{proposition}

{\bf Proof.} $i)\Rightarrow ii)$ This is Definition and Remark 2.1
vi).

$ii)\Rightarrow  iii)$ It is clear.

$iii)\Rightarrow iv)$ By [{\bf 18}, Proposition 10.1], for every
$\fp\in \supp_R(M)$, we have
$$\depth_R(M)\leq \depth_{R_{\fp}}(M_{\fp})+\dim
(R/\fp) \leq \Ht(\fp)+\dim(R/\fp)\leq \dim(R).$$ Therefore,
$\depth_{R_{\fp}}(M_{\fp})+\dim(R/\fp)=\dim(R)$, for all $\fp\in
\supp_R(M).$ On the other hand, since $M$ is complete, we have
$\fm M\neq M.$ Now, by [{\bf 20}, Theorem 3.3], we deduce that $M$
is b.b.C.M.

$iv)\Rightarrow v)$ This follows from [{\bf 16}, Corollary 3.7].

$v)\Rightarrow vi)$ This follows from [{\bf 17}, Theorem 3.6].

$vi)\Rightarrow vii)$ $M$ is b.b.C.M. by [{\bf 16}, Corollary
3.7]. Hence, by Definition and Remark 2.1 vi), we have
$\Ht_M({\fm})=\dim(R)$.

$vii)\Rightarrow iii)$ Since $M$ is a complete $R$-module, $\fm M
\neq M$. Therefore by Definition and Remark 2.1 ii), $\depth_R(M)$
is finite, So $\fm \in \supp_R(M)$. Hence, by [{\bf 17}, Corollary
3.5], we have $\depth_R(M)=\Ht_M({\fm})$.

$iv)\Rightarrow i)$ It is clear. $\Box$

\begin{lemma} Let $(R,\fm)$ be a Cohen-Macaulay local ring and $M$ a
b.b.C.M. $R$-module of finite injective dimension. Then, for every
$i>\dim_R(M)$, we have $\mu^ i_R(\fm,M)=0$.
 \end{lemma}

{\bf Proof.} We prove the lemma by induction on $\depth(R)$. Let
$\depth(R)=0$. Then $\fm \in \Ass(R)$. Thus we get an exact
sequence $$0\lo R/{\fm}\lo R\lo C\lo 0,$$ which implies the
following exact sequence $$\Ext^ i_R(R,M)\lo \Ext^
i_R(R/{\fm},M)\lo \Ext^{i+1}_R(C,M),$$ for each $i\geq 0$. Let
$r=\injdim_R(M)$. If $r> 0$, then $\Ext^ r_R(R,M)=0$, and
consequently $\Ext^ r_R(R/{\fm},M)=0$. Since $\dim(R)=0$, $\fm$ is
the unique prime ideal of $R$. Hence $\injdim_R(M)<r$ by [{\bf 2},
Corollary 3.1.12], which is a contradiction. The assertion follows
from $r=0$. Now suppose, inductively, that $\depth(R)>0$, and the
result has been proved for the smaller value. Since $\depth(R)>0$,
there exists $x \in \fm$ which is an $R$-sequence and
$M$-sequence. Since $M$ is a b.b.C.M. $R$-module, $M/xM$ is a
b.b.C.M. $R/xR$-module by [{\bf 15}, Lemma 2.3]. Therefore, in
view of Definition and Remark 2.1 vi), we have
$$\dim_{R/xR}(M/xM)=\dim(R/xR)=\dim(R)-1=\dim_R(M)-1.$$
Also, it follows from [{\bf 9}, p. 140, Lemma 2(i)] that
$$\Ext^ i_{R/xR}(R/{\fm},M/xM)\cong \Ext^{i+1}_R(R/{\fm},M),$$ for
all $i\geq 0$ and $\injdim_{R/xR}(M/xM)<\infty$. Thus, by
inductive hypothesis, for $i> \dim_R(M)-1=\dim_{R/xR}(M/xM)$, we
have $\Ext^{i+1}_R(R/{\fm},M)=0$. The assertion follows from this. $\Box$ \\

We are now ready to state and prove the main result of this
section.

\begin{theorem} Let $(R,\fm)$ be a Noetherian local ring and $M$ a
nonzero complete $R$-module. Let $\mathcal{D}(R)$ (resp.
$\mathcal{H}(M)$) be the dimension filtration of $\Spec(R)$ (resp.
$M$-height filtration of $\Spec(R)$). Then the following
conditions are equivalent:\\
i) $C(\mathcal{D}(R), M)$ is a minimal injective resolution for
$M$.\\
ii) $C(\mathcal{H}(M), M)$ is a minimal injective resolution for
$M$ and $\Ht_M(\fm)=\dim(R)$.\\
iii) $\Ht_M(\fm)=\dim(R)$ and for every $\fp \in \supp_R(M)$,
$\mu_R^ i(\fp, M)\neq 0$ if and only if $i\neq
\Ht_{M_{\fp}}({\fp}R_{\fp})$.\\
iv) $\mu_R^ i(\fm, M)=0$ if and only if $i\neq \dim(R)$.\\
v) $M$ is a b.b.C.M. $R$-module of finite injective dimension.\\
vi) $M$ is a b.C.M. $R$-module with respect to a system of
parameters for $R$ and $M$ has finite injective dimension.\\
vii) $\depth_R(M)=\dim(R)$ and $M$ has finite injective
dimension.\\
viii)
$\injdim_R(M)=\depth_R(M)=\Ht_{M}(\fm)=\dim_R(M)=\dim(R)=\depth(R)$.
\end{theorem}

{\bf Proof.} $i)\Leftrightarrow ii)$ This follows from Proposition
3.2.

$ii)\Rightarrow iii)$ Let $\fp \in \supp_{R}(M)$. We denote by
$\mathcal{H}(M)_{\fp}$ the filtration $(F_i)_{i\in \mathbb{N}_0}$
of $\Spec(R_{\fp})$, where $F_i=\{\ \fq R_{\fp}\mid \fq \in H_i \
and \ \fq \subseteq \fp \}$. Then, by [{\bf 12}, Theorem 3.5], we
have $$C(\mathcal{H}(M)_{\fp},M_{\fp})\cong
C(\mathcal{H}(M),M)_{\fp}.$$ So, $\injdim_{R_{\fp}}(M_{\fp})\leq
\Ht_{M_{\fp}}({\fp}R_{\fp})$. On the other hand, since
$C(\mathcal{H}(M),M)$ is exact, it follows from [{\bf 17},
Corollary 3.5] that
$\Ht_{M_{\fp}}({\fp}R_{\fp})=\depth_{R_{\fp}}(M_{\fp})$. Therefore
$\injdim_{R_{\fp}}(M_{\fp})=\depth_{R_{\fp}}(M_{\fp})$. This ends
this part.

$iii)\Rightarrow ii)$ Use the notation $$0\overset{b^{-2}}\lo
M\overset{b^{-1}} \lo B^{0}\lo B^{1}\lo \cdots \lo
B^{n}\overset{b^{n}} \lo B^{n+1}\lo \cdots
$$ for the Cousin complex $C(\mathcal{H}(M),M)$. By, [{\bf 17}, Corollary
3.5], $C(\mathcal{H}(M),M)$ is exact. Let $n \in \mathbb{N}_{0}$
and $(\coker b^{n-2})_{\fq}$ be the nonzero direct summand of
$B^{n}$. It is enough to show that $(\coker b^{n-2})_{\fq}$ is
injective as an $R_{\fq}$-module. Note that every injective
$R_{\fq}$-module is an injective $R$-module.

We have $\Ht_M(\fq)=n$ and $\fq \in \Supp_{R}M$. Hence, by [{\bf
17}, Corollary 2.6], $(\coker b^{n-2})_{\fq}$ is
$R_{\fq}$-isomorphic to $\COH^{n}_{\fq R_{\fq}}(M_{\fq})$, the
n-th local cohomology module of $M_{\fq}$ with respect to the
maximal ideal of the $R_{\fq}$. Since $\COH^{n}_{\fq
R_{\fq}}(M_{\fq})\neq 0$, ${\fq R_{\fq}}M_{\fq}\neq M_{\fq}$, and
so $\fq \in \supp_{R}M$. The minimal injective resolution $$0\lo
M_{\fq}\lo E^{0}\lo E^{1}\lo \cdots,$$ induces the following
complex  $$0\lo \Gamma_{\fq R_{\fq}}(E^{0})\lo \Gamma_{\fq
R_{\fq}}(E^{1})\lo \cdots \lo \Gamma_{\fq R_{\fq}}(E^{n})\lo
\Gamma_{\fq R_{\fq}}(E^{n+1})\lo \cdots .$$ By the assumption, we
have $\mu_{R}^{i}(\fq R_{\fq}, M_{\fq})=0$ for $i\neq n$. We
obtain $\Gamma_{\fq R_{\fq}}(E^{i})=0$ for every $i\neq n$.
Therefore, we have
$$\Gamma_{\fq R_{\fq}}(E^{n})\cong \COH^{n}_{\fq R_{\fq}}(M_{\fq}).$$
Hence, $\COH^{n}_{\fq R_{\fq}}(M_{\fq})$ is injective as an
$R_{\fq}$-module.

$iii)\Rightarrow iv)$ The claim follows from the fact that $\fm
\in \supp_R(M)$ (see the proof of Proposition 3.2 $vii)\Rightarrow
iii)$).

$iv)\Rightarrow v)$ By assumption and Remark 2.3 iii),
$\depth_R(M)=\dim(R)=\injdim_R(M)$. So, the claim follows from
Proposition 3.2.

$v)\Rightarrow vi)$ This is clear.

$vi)\Rightarrow vii)$ This follows from Definition and Remark 2.1
vi).

$vii)\Rightarrow viii)$ By Proposition 3.2,
$\depth_R(M)=\Ht_{M}(\fm)=\dim_R(M)=\dim(R)$. Also, we have
$\depth_R(M)\leq \injdim_R(M)$. So
$$\depth_R(M)=\Ht_{M}(\fm)=\dim_R(M)=\dim(R)\leq \injdim_R(M).$$ The assertion
follows from Remark 2.3 iv).

$viii)\Rightarrow iii)$ Proposition 3.2 yields that $M$ is
b.b.C.M.. Since  $\depth(R)=\dim(R)$, $R$ is Cohen-Macaulay and
therefore $R$ is catenary. Therefore by [{\bf 21}, Proposition
2.6], $M_{\fp}$ is b.b.C.M. as an $R_{\fp}$-module. So, from
Definition and Remark 2.1 vi),
$\depth_{R_{\fp}}(M_{\fp})=\Ht_{M_{\fp}}({\fp}R_{\fp})=\dim_{R_{\fp}}(M_{\fp})=\dim(R_{\fp})$.
On the other hand we know that
$\injdim_{R_{\fp}}(M_{\fp})<\infty$. So by Lemma 3.3, $\mu^
i_{R_{\fp}}({\fp}R_{\fp},M_{\fp})=0$ for all
$i>\dim_{R_{\fp}}(M_{\fp})$. Therefore $\mu_R^ i(\fp, M)\neq 0$ if
and only if $i\neq \Ht_{M_{\fp}}({\fp}R_{\fp})$. Note that $\mu^
i_{R_{\fp}}({\fp}R_{\fp},M_{\fp})=\mu^ i_{R}(\fp,M)$ for all $i$.
$\Box$

\section{Some characterizations for finite type complete big Cohen-Macaulay modules
of finite injective dimension}

In this section, we present some characterizations for a
Gorenstein $\hat{R}$-module $M$.

\begin{lemma} Let $(R,\fm)$ be a Noetherian local ring and $M$
an $R$-module. Suppose that $x\in \fm$ is both $R$-regular and
$M$-regular. Then:\\
i) If $M$ is a b.b.C.M. $R$-module and $\injdim_R(M)=\depth(R)$,
then
$$\injdim_{R/xR}(M/xM)=\injdim_R(M)-1.$$ ii)
$\mu_R^i(\fm,M)=\mu_{R/xR}^{i-1}(\fm/xR,M/xM)$ for all $i\geq 1$.
\end{lemma}
{\bf Proof.} i) By [{\bf 19}, Corollary 3.3.6(i)], we have
$\injdim_{R/xR}(M/xM)\leq \injdim_R(M)-1=\depth(R)-1$.

It is clear that $M/xM$ is a b.b.C.M. as an $R/xR$-module. So, by
Definition and Remark 2.1 vi), we have
$$\depth(R)-1\leq \dim(R/xR)=\depth_{R/xR}(M/xM)\leq \injdim_{R/xR}(M/xM).$$
We can deduce the result from the above facts.

ii) This follows from [{\bf 9}, p. 140, Lemma 2(i)]. $\Box$ \\

Now we are able to prove our main result of this section.

\begin{theorem} Let $(R,\fm)$ be a Noetherian local ring of dimension $d$ and $M$ a
nonzero $R$-module. The following conditions are
equivalent:\\
i) $M$ is complete and big Cohen-Macaulay with respect to a system
of parameters $x_1,\ldots,x_d$ such that $\injdim_RM$ and $\mu_R^
d(\fm,M)$ are finite.\\
ii) $M$ is a Gorenstein $\hat{R}$-module.\\
iii) $M$ is complete and for all Cohen-Macaulay $R$-module $L$ of
dimension $t$, we have

\ \ \ a) $\Ext_R^ i(L,M)=0$ for all $i\neq d-t$, and

\ \ \ b) $\Ext_R^{d-t}(\Ext_R^{d-t}(L,M),M)\cong L\otimes_RF$,
where $F$ is a finitely generated free $\hat{R}$-module.
\end{theorem}

{\bf Proof.} $i)\Rightarrow ii)$ Since $M$ is a complete
$R$-module and $\hat{M}$ is an $\hat{R}$-module, $M$ has an
$\hat{R}$-module structure. From Remark 2.3 iv) and v),
$\injdim_R(M)=\depth(R)$ and $R$ is a Cohen-Macaulay local ring.
Also, Proposition 3.2 implies that $M$ is a b.b.C.M. $R$-module.
Hence, by Lemma 4.1,
$$\mu_R^ d(\fm, M)=\mu_{R/{(x_1,\ldots,x_d)R}}^ 0({\fm}/{(x_1,\ldots,x_d)},
M/{(x_1,\ldots,x_d)M})$$ and
$$\injdim_{R/{(x_1,\ldots,x_d)R}}(M/{(x_1,\ldots, x_d)M})=
\injdim_R(M)-d=0.$$ Therefore, by Matlis decomposition theorem, we
get the following $R/{(x_1,\ldots,x_d)}R$-isomorphism.
$$M/{(x_1,\ldots,x_d)M}\cong
(E_{R/{(x_1,\ldots,x_d)R}}(R/{\fm}))^{\mu_R^ d(\fm, M)}.$$ Since
$E_{R/{(x_1,\ldots,x_d)R}}(R/{\fm})$ is finitely generated as an
$R/{(x_1,\ldots,x_d)R}$-module, it turns out that
$M/{(x_1,\ldots,x_d)M}$ is a finitely generated
$R/{(x_1,\ldots,x_d)R}$-module. Now, from the natural isomorphism
$$R/{(x_1,\ldots,x_d)R}\cong
{\hat{R}}/{(x_1,\ldots,x_d){\hat{R}}},$$ we conclude that
$M/{(x_1,\ldots,x_d)M}$ is a finitely generated
${\hat{R}}/{(x_1,\ldots,x_d){\hat{R}}}$-module. Thus, there exist
$m_1, \ldots, m_t$ in $M$ such that $M=\Sigma_{i=1}^{t}
\hat{R}m_i+\fm M$. It is easy to see that $M$ is complete as an
$\hat{R}$-module with respect to $\fm \hat{R}$-adic topology.
Also, the finitely generated $\hat{R}$-module $\Sigma_{i=1}^{t}
\hat{R}m_i$ is complete with respect to $\fm \hat{R}$-adic
topology. So, by Remark 2.3 i), we have
$M=\Sigma_{i=1}^{t}\hat{R}m_i$, that is, $M$ is a finitely
generated $\hat{R}$-module.

Let $0\lo M\lo E^{\bullet}$ be an injective resolution of $M$ as
an $\hat{R}$-module. Since any injective $\hat{R}$-module is
injective as an $R$-module, we have
\begin{equation*}
\begin{split}
\Ext_R^{i}(R/{\fm},M)&\cong
H^{i}(\Hom_{R}(R/{\fm},E^{\bullet}))\\&\cong
H^{i}(\Hom_{R}(R/{\fm},\Hom_{\hat{R}}(\hat{R},E^{\bullet})))\\&\cong
H^{i}(\Hom_{\hat{R}}({\hat{R}}/{{\fm}\hat{R}},E^{\bullet}))\\&\cong
\Ext_{\hat{R}}^{i}({\hat{R}}/{\fm}{\hat{R}},M).
\end{split}
\end{equation*}
Hence $\mu_R^ i(\fm, M)=\mu_{\hat{R}}^ i(\fm \hat{R}, M)$ for all
$i\geq 0$. On the other hand, by Theorem 3.4, we have $\mu_R^
i(\fm, M)=0$ for all $i\neq d$ and $\mu_R^ d(\fm, M)\neq 0$. So
the claim follows from [{\bf 13}, Theorem 3.11].

$ii)\Rightarrow iii)$ It is easy to see that $M$ is a complete
$R$-module with respect to $\fm$-adic topology. Since $M$ is a
Gorenstein $\hat{R}$-module, $\hat{R}$ is a Cohen-Macaulay ring by
[{\bf 13}, Theorem 3.11]. So $R$ is a Cohen-Macaulay ring.

Let $K$ be the canonical module of $\hat{R}$. As
$\hat{R}$-modules, $M$ is isomorphic to $
\oplus_{\mu_{\hat{R}}^d(\fm{\hat{R}},M)}K$ by [{\bf 14}, Corollary
2.7]. So, the claim (a) follows from Corollary 2.11 ii).

We obtain
\begin{equation*}
\begin{split}
\Ext_R^{d-t}(\Ext_R^{d-t}(L,M),M)&\cong
\Ext_R^{d-t}(\oplus_{\mu_{\hat{R}}^d(\fm{\hat{R}},M)}{\Ext_R^{d-t}(L,K),M)}\\&\cong
\oplus_{{{\mu_{\hat{R}}^d(\fm{\hat{R}},M)}}}\Ext_R^{d-t}({\Ext_R^{d-t}(L,K),M)}\\&\cong
\oplus_{{{(\mu_{\hat{R}}^d(\fm{\hat{R}},M))}}^2}\Ext_R^{d-t}({\Ext_R^{d-t}(L,K),K)}.
\end{split}
\end{equation*}
From Corollary 2.11 iii), we have
\begin{equation*}
\begin{split}
\Ext_R^{d-t}(\Ext_R^{d-t}(L,M),M)&\cong
\oplus_{{{(\mu_{\hat{R}}^d(\fm{\hat{R}},M))}}^2}(L\otimes_{R}{\hat{R}})\\&\cong
L\otimes_{R}(\oplus_{{{(\mu_{\hat{R}}^d(\fm{\hat{R}},M))}}^2}{\hat{R}}).
\end{split}
\end{equation*}
This ends the proof of b).

$iii)\Rightarrow i)$ Let $L=R/{\fm}$. Since $\dim_RL=0$, $L$ is
Cohen-Macaulay of dimension zero. Then a) and Remark 2.3 iii)
implies that $\depth_R(M)=\injdim_R(M)=d$, and so by Proposition
3.2, $M$ is b.b.C.M. On the other hand b) yields that
$$\dim_{R/{\fm}}(\Ext_R^{d}(\Ext_R^{d}(R/{\fm},M),M))$$ is finite.
We set $V=\Ext_R^{d}(R/{\fm},M)$. Then $V=\oplus_{i\in I}F_i$ such
that $F_i=R/{\fm}$ for all $i\in I$. So,
$\Ext_R^{d}(\Ext_R^{d}(R/{\fm},M),M)=\prod_{i\in I}V_i$ such that
$V_i=V$ for all $i\in I$. Thus $I$ is finite, and so
$\dim_{R/{\fm}}V$ is finite. The result follows as required.
$\Box$


\end{document}